\documentclass{pja00}

\usepackage{amsthm, amssymb, amsfonts, mathrsfs, amsmath}
\usepackage[T5]{fontenc}

\usepackage{listings}

\theoremstyle{definition}
\newtheorem{defn}{Definition}[section]
\newtheorem{rema}[defn]{Remark}
\newtheorem{nota}[defn]{Notation}

\newtheorem*{acknow}{Acknowledgments}

\theoremstyle{plain}
\newtheorem{thm}[defn]{Theorem}

\newtheorem{conj}[defn]{Conjecture}

\newtheorem{corl}[defn]{Corollary}

 \runninghead{}

\usepackage[x11names]{xcolor}
\usepackage[
  colorlinks,
]{hyperref}

\AtBeginDocument{\hypersetup{citecolor=magenta, urlcolor=magenta, linkcolor=blue}}

\title{On the dimension of the "cohits" space $\mathbb Z_2\otimes_{\mathcal A_2} H^{*}((\mathbb RP(\infty))^{\times t}, \mathbb Z_2)$\\ and some applications}

\Author{\DJ\d \abreve ng}{V\~o Ph\'uc}

\affiliation{1}{Faculty of Education Studies, University of Khanh Hoa, Vietnam}


\begin{document}

\maketitle
\begin{abstract}
We denote by $\mathbb Z_2$ the prime field of two elements and by $P_t = \mathbb Z_2[x_1, \ldots, x_t]$ the polynomial algebra of $t$ generators $x_1, \ldots, x_t$ with the degree of each $x_i$ being one. Let $\mathcal A_2$ be the Steenrod algebra over $\mathbb Z_2.$ A central problem of homotopy theory is to determine a minimal set of generators for the "cohits" space $\mathbb Z_2\otimes_{\mathcal A_2} P_t.$ This problem, which is called \textit{the "hit" problem} for Steenrod algebra, has been systematically studied for $t\leq 4.$ The present paper is devoted to the investigation of the structure of $\mathbb Z_2\otimes_{\mathcal A_2} P_t$ in some certain "generic" degrees. More specifically, we explicitly determine a monomial basis of $\mathbb Z_2\otimes_{\mathcal A_2} P_5$ in degree \mbox{$n_s=5(2^{s}-1) + 42.2^{s}$} for every non-negative integer $s.$ As a result, it confirms Sum's conjecture \cite{N.S2} for a relation between the minimal sets of $\mathcal A_2$-generators of the algebras $P_{t-1}$ and $P_{t}$ in the case $t=5$ and degree $n_s$. 
As applications, we obtain the dimension of $\mathbb Z_2\otimes_{\mathcal A_2} P_6$ in the generic degree $5(2^{s+5}-1) + n_0.2^{s+5}$ for all $s\geq 0,$ and show that the Singer's cohomological transfer \cite{W.S1} is an isomorphism in bidegree $(5, 5+n_s)$.
\end{abstract}

\section{Introduction}\label{s1} 

Recall that if $R$ is a commutative ring, then there is an $E_{\infty}$ -ring spectrum $HR$ which represents cohomology with  coeffcients in $R.$ Consider $R = \mathbb Z_2,$ the prime field of two elements, then as well known, $H^{*}(H\mathbb Z_2, \mathbb Z_2) = \mathcal A_2$ is the Steenrod algebra of $\mathbb Z_2$-cohomology operations. This is the graded homotopy of the $H\mathbb Z_2$ -module Map($H\mathbb Z_2 ,H\mathbb Z_2$), with (noncommutative) multiplication coming from a composition of cohomology operations. The algebra $\mathcal A_2$ is known to be a Hopf algebra (i.e., it is a bialgebra over the filed $\mathbb Z_2$ that has an associative multiplication and a coassociative comultiplication and is equipped with a counit, a unit and an antipode).  As an algebra over $\mathbb Z_2,$ $\mathcal A_2$  is  generated by the Steenrod squaring operations $Sq^{i}$ of degree $i$ for $i\geq 0.$ These squares act on $H^{*}(X, \mathbb Z_2)$ for any spectrum $X.$ 	Let $X= \mathbb RP(\infty)$ be the infinite real projective space. Then, it has been shown that $H^{*}(\mathbb RP(\infty), \mathbb Z_2)\cong \mathbb Z_2[x]$ and that the (left) action of $\mathcal A_2$ on $\mathbb Z_2[x]$ can be described by the rule $Sq^{i}(x^{a}) = \binom{a}{i}x^{a+i}.$ This action is deduced by induction from the case $a = 1$ by using the comultiplication on $\mathcal A_2.$ Also based on the comultiplication, one can explicitly depicte the action of $\mathcal A_2$ on \mbox{$H^{*}((\mathbb RP(\infty))^{\times t}, \mathbb Z_2)\cong P_t:= \mathbb Z_2[x_1, \ldots, x_t]$}. So, as well known, $P_t$ is an unstable left $\mathcal A_2$-module. It should be noted that $\mathcal A_2$ is not an unstable module over itself in general since $Sq^{1}Sq^{0} = Sq^{1}\neq 0.$ The $\mathbb Z_2$-algebra $P_t$ is a connected $\mathbb Z$-graded algebra. That is, $P_t = \mathbb Z_2\bigoplus (\bigoplus_{n > 0}(P_t)_n),$ where $(P_t)_n$ is the vector space of homogeneous polynomials of degree $n$. It is natural to ask: What is a minimal set of generators for $P_t$ as an $\mathcal A_2$-module? When we fix the degrees $n$ and the number of variables $t$, answering that question is the same as the problem of elementary linear algebra of finding the dimension of the quotient space, \mbox{$(QP_t)_n := (\mathbb Z_2\otimes_{\mathcal A_2} P_t)_n \cong (P_t)_n/\sum_{i > 0}{\rm Im}(Sq^{i})$}. Here $\sum_{i > 0}{\rm Im}(Sq^{i}) = (P_t)_n\cap {\rm Ker}(\epsilon)P_t,$ in which $\epsilon$ denotes the counit $\mathcal A_2\to \mathbb Z_2,$ and the Steenrod squaring operations $Sq^{i}: (P_t)_{n-i}\to (P_t)_n$ are defined, and their basic properties, such as Cartan's formula, are established. People have also shown that $(QP_t)_n$ has a module structure over the group algebra $\mathbb Z_2GL_t$ of the usual general linear group $GL_t$. The above problem, which has first appeared in the work of Peterson \cite{F.P}, is usually called the \textit{hit problem} for the Steenrod algebra. It is currently famously difficult and is of great interest to many authors.  (See Kameko \cite{M.K}, Mothebe-Uys \cite{M.M}, the present author \cite{D.P2, D.P4, Phuc}, Singer \cite{W.S2}, Sum \cite{N.S1, N.S2}, Wood \cite{R.W}, and others.)

A qualitative approach to the hit problem is that one gives conditions on elements of $(P_t)_n$ to show that they go to zero in $(QP_t)_n,$ i.e., belong to $\sum_{i > 0}{\rm Im}(Sq^{i}).$ By this, let denote $\alpha(n)$ the number of $1$'s in dyadic expansion of a non-negative integer $n.$ It is important to consider the arithmetic function $\mu: \mathbb N\to \mathbb N,$ which is defined by $\mu(0) = 0$ and $\mu(n) = \mbox{min}\{\bullet\in \mathbb N:\ \alpha(n + \bullet)\leq \bullet\}.$ Then,  Peterson's conjecture \cite{F.P}, which has been established by Wood \cite{R.W}, points out that $(QP_t)_n = 0$ when $\mu(n) > t.$ We are interested in a relatively efficient tool often used in studying hit problems. It is the $\mathbb Z_2GL_t$-modules homomorphism $\widetilde {Sq^0_*}: (QP_t)_{2n+t}\to (QP_t)_{n},$ determined by $\widetilde {Sq^0_*}([f]) = g$ if $f = x_1\ldots x_tg^{2},$ and $\widetilde {Sq^0_*}([f]) = 0$ otherwise. This linear map, given by the Kameko thesis \cite{M.K}, is related to the above $\mu$-function through the following technical result.
\begin{thm}[see Kameko \cite{M.K}]\label{dlK}
If $\mu(2n+t) = t,$ then  $\widetilde {Sq^0_*}$ is an isomorphism.
\end{thm}
By straightforward calculations using this theorem and the Peterson conjecture, it follows that the structure of $(QP_t)_n$ can be determined in each degree $n$ of the so-called the "generic" form:
\begin{equation}\label{ct}
n = r(2^s-1) + m.2^s,
\end{equation}
whenever $0\leq \mu(m) < r < t,$ and $s \geq 0.$ The cases where $t\leq 4$ are known by Peterson \cite{F.P}, Kameko \cite{M.K} and Sum \cite{N.S1}. For $t = 5$, it has been partially studied by us \cite{P.S1, D.P2, D.P4, Phuc, N.S2} and others. The problem is not yet solved in general.  The present Note is motivated by the previous results and we here continue this line of investigation by considering $(QP_t)_n$ in degrees $n$ of the form \eqref{ct}. More precisely, when $t = 5,$ we explicitly determine the structure of $(QP_5)_n$ for \mbox{$n_s :=n= 5(2^{s}-1) + 42.2^{s}$}. As a result, it confirms Sum's conjecture \cite{N.S2} for the relationship between the minimal sets of $\mathcal A_2$-generators of the polynomial rings $P_{t-1}$ and $P_{t}$ in the case $t=5$ and degree $n_s.$ An efficient approach to solve the hit problem of five variables had been provided. 
As applications, we obtain the dimension of $(QP_6)_n$ in degree \mbox{$5(2^{s+5}-1) + n_0.2^{s+5}$} for all $s\geq 0$, and show that Singer's conjecture \cite{W.S1} for the fifth cohomological transfer holds in degrees $n_s.$    

This paper is an announcement whose details will appear elsewhere.

\section{The Necessary Preliminaries}

In order to formulate the main results of this text and for the convenience of the reader, let us describe some preliminary material.

\begin{defn}[{\bf Weight vector and exponent vector}] We say that a sequence of non-negative integers \mbox{$\omega = (\omega_1, \omega_2, \ldots, \omega_i,\ldots)$} is a \textit{weight vector}, if $\omega_i  = 0,$ for $i\gg 0.$ Then, one defines \mbox{$\deg(\omega) = \sum_{i\geq 1}2^{i-1}\omega_i,$} With a natural number $n,$ let denote $\alpha_j(n)$ the $j$-th coefficients in dyadic expansion of $n,$ then $\alpha(n) = \sum_{j\geq 0}\alpha_j(n),$ and $n = \sum_{j\geq 0}\alpha_j(n)2^j,$ where $\alpha_j(n)\in \{0, 1\}$ for all $j\geq 0.$ For a monomial \mbox{$x = x_1^{a_1}x_2^{a_2}\ldots x_t^{a_t}\in P_t,$} let us consider two sequences $\omega(x) :=(\omega_1(x), \ldots, \omega_i(x), \ldots)$ and $(a_1, a_2, \ldots, a_t),$ where $\omega_i(x)=\sum_{1\leq j\leq t}\alpha_{i-1}(a_j),$ for every $i.$ They are respectively called the {\it weight vector} and the \textit{exponent vector} of $x.$  By convention, the sets of all the weight vectors and the exponent vectors are given the left lexicographical order.
\end{defn}

\begin{defn}[{\bf Linear order on {\boldmath{$P_t$}}}]
Let $x = x_1^{a_1}x_2^{a_2}\ldots x_t^{a_t}$ and \mbox{$y = x_1^{b_1}x_2^{b_2}\ldots x_t^{b_t}$} be the monomials in $(P_t)_n.$ We write $a,\, b$ for the exponent vectors of $x$ and $y,$ respectively. We say that $x  < y$ if and only if one of the following holds:
\begin{enumerate}
\item[(i)] $\omega(x) < \omega(y);$
\item[(ii)] $\omega(x) = \omega(y)$ and $a < b.$
\end{enumerate}
\end{defn}

\begin{defn}[{\bf Equivalence relations on {\boldmath{$P_t$}}}]  For a weight vector $\omega$ of degree $n,$ let us denote two subspaces associated with $\omega$ by\\[1mm]
$\begin{array}{ll}
(P_t)_n^{\omega} &= \langle\{ x\in (P_t)_n|\, \deg(x) = \deg(\omega),\  \omega(x)\leq \omega\}\rangle,\\
(P_t)_n^{< \omega} &= \langle \{ x\in (P_t)_n|\, \deg(x) = \deg(\omega),\  \omega(x) < \omega\}\rangle.
\end{array}$\\[1mm]
Let $f$ and $g$ be two homogeneous polynomials in $(P_t)_n,$ we define the binary relations "$\sim$" and "$\sim_{\omega}$" on $(P_t)_n$:
\begin{enumerate}
\item [(i)] $f \sim g $ if and only if $(f - g)\in \sum_{i > 0}{\rm Im}(Sq^{i})$. (Note that since working mod 2, $f-g = f+g.$)
\item[(ii)]  $f \sim_{\omega} g$ if and only if $f, \, g\in (P_t)_n^{\omega}$ and\\[1mm] 
$(f -g)\in (((P_t)_n^{\omega}\cap \sum_{i > 0}{\rm Im}(Sq^{i}))  + (P_t)_n^{< \omega}).$\\[1mm]
In particular,  if $f\sim 0$ (resp. $f\sim_{\omega} 0$), then we say that $f$ is \textit{hit} (resp. \textit{$\omega$-hit}). It is straightforward to see that the above binary relations are equivalence ones. So, one has a quotient space $$(Q_n^{\otimes t})^{\omega} := (P_t)_n^{\omega}/ (( (P_t)_n^{\omega})\cap \sum_{i > 0}{\rm Im}(Sq^{i})) + (P_t)_n^{< \omega}),$$ from which it is not difficult to check the isomorphism  $(QP_t)_n\cong \bigoplus_{\deg(\omega) = n}(Q_n^{\otimes t})^{\omega}$ (see \cite{D.P2}).
\end{enumerate}
\end{defn}

\begin{defn}[{\bf Admissible monomial and inadmissible monomial}] A monomial $x$ in $(P_t)_n$ is said to be {\it inadmissible} if there exist monomials $y_1, y_2,\ldots, y_k$ in $(P_t)_n$ such that $y_j < x$ for all $j,\ 1\leq j\leq k,$ and $x \sim  (y_1 + y_2 + \cdots + y_k).$ Then, we say that $x$ is {\it admissible}, if it is not inadmissible.
\end{defn}

The following important theorem, which is due to Kameko \cite{M.K}, can be viewed as a criterion in the calculation of the inadmissible monomials.

\begin{thm}\label{dlKS}
Let $x$ and $z$ be monomials in $(P_t)_n.$ For an integer $r >0,$ assume that there exists an index $i>r$ such that $\omega_i(x) = 0$, from which if $z$ is inadmissible, then $xz^{2^r}$ is, too.
\end{thm} 

From the above data, we have seen that $(QP_t)_n$ is an $\mathbb Z_2$-vector space with a basis consisting of all the classes represent by the admissible monomials in $(P_t)_n.$

\begin{defn}[{\bf Spike monomial}]
A monomial $z = \prod_{1\leq j\leq t}x_j^{b_j}$ in $(P_t)_n$ is called a {\it spike} if every exponent $b_j$ is of the form $2^{\beta_j} - 1.$ In particular, if the exponents $\beta_j$ can be arranged to satisfy $\beta_1 > \beta_2 > \ldots > \beta_{r-1}\geq \beta_r \geq 1,$ where only the last two smallest exponents can be equal, and $\beta_j = 0$ for $ r < j  \leq t,$ then $z$ is called a {\it minimal spike}.

\end{defn}

\begin{thm}[see Ph\'uc-Sum \cite{P.S1}]\label{dlPS}
All the spikes in $(P_t)_n$ are admissible and their weight vectors are weakly decreasing. Furthermore, if a weight vector $\omega = (\omega_1, \omega_2, \ldots)$ is weakly decreasing and \mbox{$\omega_1\leq t,$} then there is a spike $z\in (P_t)_n$ such that $\omega(z) = \omega.$
\end{thm}

\begin{thm}[see Singer \cite{W.S2}]\label{dlSi}
Suppose that $x\in (P_t)_n$ and $\mu(n)\leq t.$ Let $z$ be the minimal spike in $(P_t)_n.$ If $\omega(x) < \omega(z),$ then $x\sim 0.$
\end{thm}

For $1\leq l\leq t,$ let $\mathsf{q}_{(l,\,t)}: P_{t-1}\to P_t$ be an $\mathcal A_2$-homomorphism, which is depicted by 
$$ \mathsf{q}_{(l,\,t)}(x_j) = \left\{ \begin{array}{ll}
{x_j}&\text{if }\;1\leq j \leq l-1, \\
x_{j+1}& \text{if}\; l\leq j \leq t-1.
\end{array} \right.$$

We assume throughout this text that $\mathscr {C}^{\otimes t}_n$ is the set of all admissible monomials in $(P_t)_n.$ The following technicality is very useful for our calculations in the next section.

\begin{thm}[see Mothebe-Uys \cite{M.M}]\label{dlM}
Let $l, d$ be positive integers such that $1\leq l\leq t.$ If \mbox{$x\in \mathscr {C}^{\otimes (t-1)}_{n},$} then $x_l^{2^{d}-1}\mathsf{q}_{(l,\,t)}(x)\in \mathscr {C}^{\otimes t}_{n + 2^{d}-1}.$
\end{thm}

We set $ \mathcal N_t := \{(l, \mathscr L)\;|\; \mathscr L = (l_1,l_2,\ldots, l_r)\},$ where $1\leq l < l_1< l_2 < \ldots < l_r\leq t$ and \mbox{$0\leq r \leq t-1.$} By convention, $\mathscr L = \emptyset$ if $r = 0.$ Let $r = \ell(\mathscr L)$ be the length of $\mathscr L$. For each the pair $(l, \mathscr L)\in\mathcal{N}_t,$ and $1\leq u < r,$ let us denote by $x_{(\mathscr L,\,u)} := x_{l_u}^{2^{r-1} + 2^{r-2} +\, \cdots\, + 2^{r-u}}\prod_{u < d\leq r}x_{l_d}^{2^{r-d}},$ where $x_{(\emptyset, 1)} = 1.$ Sum \cite{N.S1} defines an $\mathbb Z_2$-linear function $\psi_{(l, \mathscr L)}: P_{t-1}\longrightarrow P_t,$ which is determined by  
$$ \begin{array}{ll}
& \psi_{(l, \mathscr L)}(\prod_{1\leq j\leq t-1}x_j^{a_j}) \\
&= (x_l^{2^{r} - 1}\mathsf{q}_{(l,\,t)}(\prod_{1\leq s\leq t-1}x_s^{a_s}))/x_{(\mathscr L,\,u)},
\end{array}$$ if there exists $u$ such that:
 \begin{equation}\label{ct2}
\begin{array}{ll}
a_{l_1 - 1} +1= \cdots = a_{l_{(u-1)} - 1} +1 = 2^{r},\\
 a_{l_{u} - 1} + 1 > 2^{r},\\
\alpha_{r-d}(a_{l_{u} - 1}) -1 = 0,\, 1\leq d\leq u, \\
\alpha_{r-d}(a_{l_{d}-1}) -1 = 0,\, \ u+1 \leq d \leq r,
\end{array}
\end{equation}
and $\psi_{(l, \mathscr L)}(\prod_{1\leq j\leq t-1}x_j^{a_j}) = 0$ otherwise. One has the following observation.

\begin{rema}
If  $\mathscr L = \emptyset,$ then $\psi_{(l, \mathscr L)} = \mathsf{q}_{(l,\,t)}$ for $1\leq l\leq t.$ It is in fact not hard to show that if $\psi_{(l, \mathscr L)}( \prod_{1\leq s\leq t-1}x_s^{a_s})\neq  0,$ then $\omega(\psi_{(l, \mathscr L)}( \prod_{1\leq s\leq t-1}x_s^{a_s})) = \omega( \prod_{1\leq s\leq t-1}x_s^{a_s}).$ Alternatively, $\psi_{(l, \mathscr L)}$ is not the homomorphism of $\mathcal A_2$-modules in general. As an illustrated example,  we consider $t = 4,$  $\mathscr L = (2,3,4)\neq \emptyset,$ and the monomial $x_1^{12}x_2^{6}x_3^{9}\in (P_3)_{27},$ then, straightforward calculations show that\\[1mm]
$\begin{array}{ll}
\psi_{(1, \mathscr L)}(x_1^{12}x_2^{6}x_3^{9}) &= \dfrac{x_1^{2^{3}-1}\mathsf{q}_{(1, 4)}(x_1^{12}x_2^{6}x_3^{9})}{x_{(\mathscr L,\,1)}} \\
&= \dfrac{x_1^{7}x_2^{12}x_3^{6}x_4^{9}}{x_2^{4}x_3^{2}x_4} = x_1^{7}x_2^{8}x_3^{4}x_4^{8}\in (P_4)_{27},
\end{array}$\\[1mm]
and $\omega(x_1^{12}x_2^{6}x_3^{9}) =  (1,1,2,2) =\omega(\psi_{(1, \mathscr L)}(x_1^{12}x_2^{6}x_3^{9})).$ So, using the Steenrod squares $Sq^{2},$ one gets  \\[1mm]
$\begin{array}{ll}
&Sq^{2}(\psi_{(1, \mathscr L)}(x_1^{12}x_2^{6}x_3^{9})) = x_1^{9}x_2^{8}x_3^{4}x_4^{8}\\
&\quad\quad\quad\quad\quad\quad \neq x_1^{7}x_2^{8}x_3^{6}x_4^{8}= \psi_{(1, \mathscr L)}(Sq^{2}(x_1^{12}x_2^{6}x_3^{9})).
\end{array}$\\[1mm]
\end{rema}
For a subset $\mathscr C\subset (P_{t-1})_n,$ let us consider the sets:\\[1mm]
$\begin{array}{ll}
\widetilde {\Phi^0}(\mathscr C) &= \bigcup\limits_{1\leq l \leq t}\psi_{(l, \emptyset)}(\mathscr  C) =  \bigcup_{1\leq l \leq t}\mathsf{q}_{(l, t)}(\mathscr C),\\
 \widetilde{\Phi^{>0}}(\mathscr  C) &= \bigcup\limits_{(l; \mathscr L)\in\mathcal{N}_t,\;1 \leq \ell(\mathscr L) \leq t-1}(\psi_{(l, \mathscr L)}(\mathscr C)\setminus (P_t)_n^{0}),
\end{array}$\\[1mm]
where $$(P_t)_n^{0} = \langle  \{\prod_{1\leq i\leq t}x_i^{a_i}\in (P_t)_n:\, \prod_{1\leq i\leq t}a_i = 0\}\rangle$$  is an $\mathbb Z_2$-subspace of $(P_t)_n.$  For later convenience, let us denote by $\widetilde{\Phi_*}(\mathscr C) := \widetilde {\Phi^0}(\mathscr C) \bigcup \widetilde {\Phi^{>0}}(\mathscr C).$ Let $\omega$ be a weight vector of degree $n$ and let \mbox{$(\mathscr{C}^{\otimes t}_{n})^{\omega}:= \mathscr {C}^{\otimes t}_{n}\cap (P_t)_{n}^{\omega}.$} In \cite{N.S2}, Sum sets up the following conjecture, which plays a crucial role in studying the minimal set of $\mathcal A_2$-generators for $P_t$ in each positive degree.

\begin{conj}\label{gtS}
Under the notations chosen, for each $(l, \mathscr L)\in \mathcal N_t,$ and \mbox{$1\leq r = \ell(\mathscr L)\leq t-1,$} if $x = \prod_{1\leq j\leq t-1}x_j^{a_j}\in \mathscr C^{\otimes (t-1)}_n,$ and there exist \mbox{$u,\ 1\leq u\leq r,$} which satisfies \eqref{ct2}, then \mbox{$ \psi_{(l, \mathscr L)}(x) \in \mathscr C^{\otimes t}_n.$} Moreover, if  $\omega$ is a weight vector of degree $n$, then $\widetilde{\Phi_*}((\mathscr{C}^{\otimes (t-1)}_{n})^{\omega})\subseteq  (\mathscr{C}^{\otimes t}_{n})^{\omega}.$ 
\end{conj}

Let us consider the set $$\mathcal U(t, n) := \{\rho\in\mathbb N\ :\ \alpha(n - (2^\rho - 1) + t-1)\leq t-1\}.$$ For each $\rho\in \mathcal U(t, n),$ suppose that the monomial $y = \prod_{1\leq i\leq t-1}x_i^{a_i}\in \mathscr C^{\otimes (t-1)}_{n-(2^\rho - 1)}.$  We denote\\[1mm] $\mathscr C(t, n):= \bigcup\limits_{\rho\in \mathcal U(t, n)}\langle\{x_1^{a_1}\ldots x_{j-1}^{a_{j-1}}x_j^{2^{\rho} - 1}x_{j+1}^{a_{j+1}}\ldots x_t^{a_{t-1}}\in (P_t)_n\ :\ 1\leq j\leq t\}\rangle$ and write $$(\mathscr C(t, n))^{\omega} = \{x\in \mathscr C(t, n):\, \omega(x) = \omega\}.$$
Combining Conjecture \ref{gtS} with a result in \cite{M.M}, we propose the following.
\begin{conj}\label{gtP}
For each $(l, \mathscr L)\in \mathcal N_t,$ if $x\in (\mathscr C_n^{\otimes (t-1)})^{\omega}$ and $$\mbox{max}\bigg\{r = \ell(\mathscr L) \leq t-1:  \psi_{(l, \mathscr L)}(x)\neq 0\bigg\} = |\mathcal U(t, n)|,$$ then
$$(\widetilde{\Phi_*}((\mathscr{C}^{\otimes (t-1)}_{n})^{\omega})\cup (\mathscr C(t, n))^{\omega})= (\mathscr{C}^{\otimes t}_{n})^{\omega}.$$ 
\end{conj}

To close this section, we recall a result of Sum \cite{N.S1} on the inductive formula for the dimension of $(QP_t)_n.$

\begin{thm}[see Sum \cite{N.S1}]\label{dlS}
Consider the degree $n$ of the form \eqref{ct} with $r = t-1,$ and $s, m$ positive integers such that $1\leq t-3\leq \mu(m)\leq t-2,$ and $\mu(m) = \alpha(m + \mu(m)).$ Then for each $s\geq t-1,$ we have $\dim (QP_t)_{n}  = (2^{t}-1)\dim (QP_{t-1})_{m},$
\end{thm}

By using iteration of  the Kameko homomorphism $\widetilde {Sq^0_*}$ and the proof of Theorem \ref{dlS}, one has an equivalent statement of Theorem \ref{dlS} that will be applied in the sequel.

\begin{thm}\label{dlP}
Let us consider the generic degrees of the form \eqref{ct}: $n_s := (t-1)(2^{s}-1) + m.2^{s}.$ Suppose there is a integer $\zeta$ such that $0\leq \zeta < s$ and $1\leq t-3\leq \mu(n_{\zeta}) = \alpha(n_{\zeta}+ \mu(n_{\zeta}))\leq t-2.$ Then, for each $s \geq \zeta,$ we have
$$ \begin{array}{ll}
&\dim (QP_t)_{(t-1)(2^{s-\zeta + t-1}-1) + n_{\zeta}2^{s-\zeta + t-1}} \\
&=(2^{t}-1)\dim (QP_{t-1})_{n_s}.
\end{array}$$
\end{thm}

\section{Statement of results}
We now survey our main results in this Note. Let us consider the generic degree of the form \eqref{ct} with $r = t = 5,$ and $m = 42.$ Then, we have $n_s:= 5(2^{s}-1) + 42.2^{s}$ with $s\geq 0.$ To solve the hit problem of five variables in this degree, our main idea is to combine the homomorphisms $\psi_{(l, \mathscr L)}$ and Theorem \ref{dlM} above. This approach, which is quite effective, helps us to reduce many computations. This is shown specifically in Theorem \ref{dlc1} below.

\begin{nota} 
In what follows:\\[1mm]
$\begin{array}{ll}
(P_t)_{n_s}^{>0} &= \langle  \{\prod_{1\leq i\leq t}x_i^{a_i}\in (P_t)_{n_s}:\, \prod_{1\leq i\leq t}a_i >0\}\rangle\\[0.7mm]
(QP_t)_{n_s}^0 &= \mathbb Z_2\otimes_{\mathcal A_2} (P_t)_{n_s}^{0}, (QP_t)_{n_s}^{>0}= \mathbb Z_2\otimes_{\mathcal A_2} (P_t)_{n_s}^{>0}.
\end{array}$\\[1mm]
Clearly, $(P_t)_{n_s}^{>0} $ is an $\mathbb Z_2$-subspace of $(P_t)_n,$ and $(QP_t)_{n_s} \cong (QP_t)_{n_s}^0\,\bigoplus\, (QP_t)_{n_s}^{>0}.$ 

For a polynomial $f\in (P_t)_{n_s},$ we denote by $[f]$ the classes in $(QP_t)_{n_s}$ represented by $f.$ If $\omega$ is a weight vector of degree $n_s,$ and $f\in (P_t)_{n_s}^{\omega},$ then denote by $[f]_\omega$ the classes in $(QP_t)_{n_s}^{\omega}$ represented by $f.$ For a subset $\mathscr{C}\subset (P_t)_{n_s},$ denote by \mbox{$[\mathscr C] = \{[f]\, :\, f\in \mathscr C\}.$} If $\mathscr C\subset  (P_t)_{n_s}^{\omega},$ then denote by $[\mathscr C]_{\omega} = \{[f]_{\omega}\, :\, f\in \mathscr C\}.$ It should be noted that if $\omega$ is a weight vector of mimimal spike, then $[\mathscr C]_{\omega} = [\mathscr C].$ We put 
$$\begin{array}{ll}
(\mathscr{C}^{\otimes t}_{n_s})^{\omega^{0}}& := (\mathscr{C}^{\otimes t}_{n_s})^{\omega}\cap  (P_t)_{n_s}^{0},\\[0.7mm]
 (\mathscr{C}^{\otimes t}_{n_s})^{\omega^{>0}} &:= (\mathscr{C}^{\otimes t}_{n_s})^{\omega}\cap (P_t)_{n_s}^{>0},\\[0.7mm]
(Q_{n_s}^{\otimes t})^{\omega^{0}}&:= (QP_t)_{n_s}^{\omega}\cap (QP_t)_{n_s}^{0},\\[0.7mm]
 (Q_{n_s}^{\otimes t})^{\omega^{>0}} &:= (QP_t)_{n_s}^{\omega}\cap (QP_t)_{n_s}^{>0}.
\end{array}$$
Then, $[(\mathscr{C}^{\otimes t}_{n_s})^{\omega}]_\omega,$ $[(\mathscr{C}^{\otimes t}_{n_s})^{\omega^{0}}]_\omega$ 
and $[(\mathscr{C}^{\otimes t}_{n_s})^{\omega^{>0}}]_\omega$ are respectively the bases of the vector spaces $(Q_{n_s}^{\otimes t})^{\omega},\ (Q_{n_s}^{\otimes t})^{\omega^{0}},$ and $(Q_{n_s}^{\otimes t})^{\omega^{>0}}.$
\end{nota}

For a start, we have the following important remark.
\begin{rema}\label{nx1}
Noticing that $\mu(n_s) = 5$ for all $s > 0.$ This, together with Theorem \ref{dlK}, implies that the map
\begin{equation}\label{pt1}
(\widetilde {Sq^0_*})^{s}: (QP_5)_{n_s}  \to (QP_5)_{n_0}
\end{equation}
is an isomorphism of $\mathbb Z_2$-vector spaces for any \mbox{$s\geq 0$}, from which we only need to compute $(QP_5)_{n_0}.$  We have $n_0 = 5(2^{0}-1) + 42.2^{0} = 42.$ Then, if \mbox{$x\in \mathscr C^{\otimes 5}_{42},$} then the weight vector $\overline{\omega}:=\omega(x)$ is one of the following sequences: $\overline{\omega}_{(1)} = (4,3,2,1,1),$ \mbox{$\overline{\omega}_{(2)}= (4,3,2,3),$} and $\overline{\omega}_{(3)} = (4,3,4,2).$ Indeed, it should be noted that \mbox{$z = x_1^{31}x_2^{7}x_3^{3}x_4\in (P_5)_{42}$} is minimal spike, and so, by Theorem \ref{dlPS}, $z$ belongs to $\mathscr C^{\otimes 5}_{42}.$ Moreover, because $\omega(z) = (4,3,2,1,1),$ based upon Theorem \ref{dlSi}, we see that if $x\in \mathscr C^{\otimes 5}_{42},$ then  $\omega_1(x) = 4.$ This means that $x$ is of the form $x_ix_jx_kx_ly^{2}$ with $y\in (P_5)_{19}$ and \mbox{$1\leq i<j<k<l\leq 5.$} Since $x\in\mathscr C^{\otimes 5}_{42},$ by Theorem \ref{dlKS}, we must have $y\in \mathscr C^{\otimes 5}_{19}.$ So, the desired conclusion now follows from a result in \cite{T1} that $\omega(y)$ belongs to the set $\{(3,2,1,1),\, (3,2,3),\, (3,4,2)\}.$ 
\end{rema}

From the above remarks, the following is immediate:
\begin{equation}\label{pt2}
 (QP_5)_{42}\cong (QP_5)_{42}^{0}\bigoplus \big(\bigoplus_{1\leq j\leq 3} (Q^{\otimes 5}_{42})^{\overline{\omega}_{(j)}^{>0}}\big).
\end{equation}
Since $\mu(42) = 4,$ $(QP_4)_{42}\cong (QP_4)_{19}.$ Then, one can verify that $$(QP_4)_{42} = (QP_4)_{42}^{>0} =  (Q_{42}^{\otimes 4})^{\overline{\omega}_{(1)}^{>0}} = (Q_{42}^{\otimes 4})^{\overline{\omega}_{(1)}},$$ and so, $\dim (QP_4)_{42} = \dim (QP_4)_{19} = 140.$ 
On the other side, it has been shown (see also \cite{D.P4}) that 
$$ \dim (QP_5)^0 = \sum_{1\leq  s\leq 4}\binom{5}{s}\dim (QP_s)_n^{>0},$$
for all $n\geq 0.$ Using this formula, together with Peterson's conjecture and Theorem \ref{dlK}, we are forced to conclude that
\begin{equation}\label{pt3}
\begin{array}{ll}
\dim (QP_5)_{42}^{0} &= \dim (Q_{42}^{\otimes 5})^{\overline{\omega}_{(1)}^{0}}\\
\medskip
& =  140.\binom{5}{4}= 700,\\
(\mathscr C_{42}^{\otimes 5})^{\overline{\omega}_{(1)}^{0}} &= \widetilde {\Phi^0}((\mathscr C^{\otimes 4}_{42})^{\overline{\omega}_{(1)}^{>0}}) = \widetilde {\Phi^0}(\mathscr C^{\otimes 4}_{42}) .
\end{array}
\end{equation}

For $1\leq l\leq 5,$ we consider the following sets:
$$ \begin{array}{ll}
\mathscr C_1 &:= \{x_l\mathsf{q}_{(l,\,5)}(x):\ x\in (\mathscr {C}^{\otimes 4}_{41})^{>0},\} ,\\[0.7mm]
\mathscr C_2 &:= \{x_l^{3}\mathsf{q}_{(l,\,5)}(x):\ x\in (\mathscr {C}^{\otimes 4}_{39})^{>0}\},\\[0.7mm]
\mathscr C_3 &:= \{x_l^{7}\mathsf{q}_{(l,\,5)}(x):\ x\in (\mathscr {C}^{\otimes 4}_{35})^{>0}\},\\[0.7mm]
\mathscr C_4 &:= \{x_l^{31}\mathsf{q}_{(l,\,5)}(x):\ x\in (\mathscr {C}^{\otimes 4}_{11})^{>0}\}.
\end{array}$$
Then, using Theorem \ref{dlM}, it can be shown that $$\mathscr C^{>0}(\overline{\omega}_{(1)}):= \bigcup_{1\leq i\leq 4}\mathscr C_i\subseteq (\mathscr C^{\otimes 5}_{42})^{\overline{\omega}_{(1)}^{>0}}.$$ We put\\[1mm]
$E:= \widetilde{\Phi^{>0}}((\mathscr C_{42}^{\otimes 4})^{\overline{\omega}_{(1)}})\setminus (\widetilde{\Phi^{>0}}((\mathscr C_{42}^{\otimes 4})^{\overline{\omega}_{(1)}})\cap \mathscr C^{>0}(\overline{\omega}_{(1)})).$ Then, we have the theorem below, which is our first main result.

\begin{thm}\label{dlc1}
The following statements hold:
\begin{itemize}
\item[i)] $(\mathscr C_{42}^{\otimes 5})^{\overline{\omega}_{(j)}^{>0}} =  (\mathscr C_{42}^{\otimes 5})^{\overline{\omega}_{(j)}}=\emptyset,$ for $j = 2, 3.$

\item[ii)]  $(\mathscr C_{42}^{\otimes 5})^{\overline{\omega}_{(1)}^{>0}} = E \bigcup F \bigcup \mathscr C^{>0}(\overline{\omega}_{(1)}),$
\end{itemize}
where $|E|= 542,\ |F|  = 248,$ and $|\mathscr C^{>0}(\overline{\omega}_{(1)})| = 1030.$ Consequently, $(Q_{42}^{\otimes 5})^{\overline{\omega}_{(1)}^{>0}}$ is $1820$-dimensional.
\end{thm}

From \eqref{pt1}, \eqref{pt3}, and Theorem \ref{dlc1},  it may be concluded that Conjecture \ref{gtS} holds for $t = 5$ and the generic degree $n_s$ for any $s\geq 0.$ 

As immediate consequences from  \eqref{pt1}, \eqref{pt2}, \eqref{pt3} and Theorem \ref{dlc1}, we may assert that

\begin{corl}\label{hq}
There exist exactly $2520$ admissible monomials of degree $n_s$ in $P_5$ for all $s\geq 0.$ Consequently, $(QP_5)_{n_s}$ is $2520$-dimensional.
\end{corl}

The above computations confirm the dimension of $(QP_5)_{n_s},$ which is informed in \cite{D.P4} by using the MAGMA computer algebra system.  

Next, based on Theorem \ref{dlP} and Corollary \ref{hq}, we have immediately the following.

\begin{corl}\label{hq3}
When degree $n_s$ is as in Theorem \ref{dlP}, for $t = 6$, $m = 42,$ and $s\geq 0,$ we have $\dim (QP_6)_{n_{s+5}}= (2^{6}-1)\dim (QP_5)_{n_s}  = 158760.$
\end{corl}

{\bf Singer's cohomological transfer.} The next contribution of the paper is to apply Theorem \ref{dlc1} into the study of the behavior of the fifth transfer homomorphism of W. Singer \cite{W.S1}. It may need to be recalled that the group $GL_t$ acts regularly on $P_t$ by matrix substitution. Further, the two actions of $\mathcal A_2$ and $GL_t$ upon $P_t$ commute with each other; hence there is an inherited action of $GL_t$ on $QP_t = \mathbb Z_2\otimes_{\mathcal A_2}P_t.$ Thus, the hit problem becomes important in studying the modular representation of $GL_t.$ This is of interest since the Singer cohomological "transfer" \cite{W.S1} relates $(QP_t)_n$ to $H^{t, t+n}(\mathcal A_2, \mathbb Z_2) = {\rm Ext}_{\mathcal A_2}^{t, t+n}(\mathbb Z_2, \mathbb Z_2),$ the $\mathbb Z_2$-cohomology groups of $\mathcal A_2,$ and thus to the stable homotopy groups of spheres. Singer's transfer is established as follows. Let us recall first that if $\mathcal M$ the category of graded left $\mathcal A_2$-modules and degree zero of $\mathcal A_2$-linear maps, then for each integer $t,$ the supension functor \mbox{$S^{t}: \mathcal M\to \mathcal M$} is defined by \mbox{$(S^{t}U)^{n} := U^{n-t}$} where \mbox{$n\in\mathbb Z$}. The action of $\mathcal A_2$ on $S^{t}U$ is given  by \mbox{$\theta(S^{t}u) = S^{t}(\theta u),$} for all $u\in U$ and $\theta\in \mathcal A_2.$ Next, consider the polynomial ring \mbox{$P_1 = \mathbb Z_2[x]$} with $\deg(x) = 1,$ the canonical $\mathcal A_2$-action on $P_1$ is extended to an $\mathcal A_2$-action on the ring of finite Laurent series $\mathbb Z_2[x_1, x_1^{-1}]$. Then, there exists an $\mathcal A_2$-submodule \mbox{$\overline{P} = \langle \{x_1^{t}|\ t\geq -1\}\rangle$} of $\mathbb Z_2[x_1, x_1^{-1}].$ One has a short-exact sequence $ 0\to P_1\xrightarrow{q} \overline{P}\xrightarrow{\pi} S^{-1}\mathbb Z_2,$ where $q$ is the inclusion and $\pi$ is given by $\pi(x_1^{t}) = 0$ if $t\neq -1$ and $\pi(x_1^{-1}) = 1.$ Writing $e_1$ for the corresponding element in ${\rm Ext}_{\mathcal A_2}^1(S^{-1}\mathbb Z_2, P_1).$ Basing the cross, the Yoneda and the cap products in co(homology) with $\mathbb Z_2$-coeffcients, we have a homomorphism
$$ \begin{array}{ll}
 Tr^{\mathcal A_2}_t: {\rm Tor}_t^{\mathcal A_2}(\mathbb Z_2, S^{-t}\mathbb Z_2) &\longrightarrow {\rm Tor}_0^{\mathcal A_2}(\mathbb Z_2, P_t) = QP_t\\
\hspace{2cm} z&\longmapsto e_t\cap z,
\end{array}$$
where $e_t = (e_1\times P_{t-1})\circ \cdots \circ (e_1\times P_1)\circ e_1$. Its image is a submodule of the space of $GL_t$-invariants $(QP_t)^{GL_t}.$ It is known, the supension $S^{-t}$ induces an isomorphism $${\rm Tor}_{t, n}^{\mathcal A_2}(\mathbb Z_2, S^{-t}\mathbb Z_2)\cong H_{t, t+n}(\mathcal A_2, \mathbb Z_2)$$ and so $Tr^{\mathcal A_2}_t$ induces the homomorphism
$$ Tr^{\mathcal A_2}_t: H_{t, t+n}(\mathcal A_2, \mathbb Z_2)\to (QP_t)_{n}^{GL_t}.$$
The dual $ (Tr^{\mathcal A_2}_t)^{*}: ((QP_t)_{n}^{GL_t})^{*}\to H^{t, t+n}(\mathcal A_2, \mathbb Z_2)$
of $Tr^{\mathcal A_2}_t$ is called \textit{the $t$-th cohomological transfer}. 
 In \cite{W.S1}, Singer proves that $ (Tr_t^{\mathcal A_2})^{*}$ is an isomorphism for $t =1,2.$ Afterwards, based on Kameko's thesis \cite{M.K}, Boardman \cite{J.B} states that $ (Tr^{\mathcal A_2}_3)^{*}$ is also an isomorphism. These works tell us that Singer's transfer map is highly non-trivial, and so, it serves as a reliable tool for depicting the cohomology of the Steenrod ring. In higher cohomological degrees, Singer \cite{W.S1} gives some calculations to show that the fourth transfer is an isomorphism up to a range of internal degrees, but the fifth transfer is not an epimorphism, from which he makes the following prediction.

\begin{conj}\label{gtSi}
The cohomological transfer is a monomorphism.
\end{conj}

This was confirmed for $1\leq t\leq 3$ by Singer himself \cite{W.S1} and Boardman \cite{J.B}. Our recent work \cite{D.P6} shows that the conjecture is also true for $t = 4.$ 
So far it remains open in general. Very little information is known when $t = 5.$  Now, based upon an admissible monomial basis for $QP_5$ in degree $5(2^{s}-1) + 42.2^{s}$ (see Theorem \ref{dlc1}), we verify Conjecture \ref{gtSi} for $t = 5$ and the respective degree. The following theorem is our second main result.

\begin{thm}\label{dlc3}
The transfer homomorphism $(Tr^{\mathcal A_2}_5)^{*}: ((QP_5)_{5(2^{s}-1) + 42.2^{s}}^{GL_5})^{*}\to H^{5, 47.2^{s}}(\mathcal A_2, \mathbb Z_2)$ is a trivial isomorphism, for all $s\geq 0.$
\end{thm}

This result can be proved by explicitly computing the space of $GL_5$-invariants $(QP_5)_{42}^{GL_5}$ combining the fact that $H^{5, 47.2^{s}}(\mathcal A_2, \mathbb Z_2)$ is trivial (see Lin \cite{W.L}). Thus, with these data in hand, the reader can see that

\begin{corl}
Conjecture \ref{gtS} holds for the rank 5 case and degree $5(2^{s}-1) + 42.2^{s}$ with $s\geq 0.$
\end{corl}

\begin{acknow}
The author research is partially funded by the NAFOSTED grant No. 101.04-2017.05. 

\end{acknow}

\end{document}